\documentclass[review]{elsarticle}
\usepackage{mathtools}
\usepackage{amsmath}
\usepackage{amssymb}
\usepackage{amsmath,esint}

\usepackage{hyperref}

\journal{Applied Mathematics and Computation}

\begin{document}

\begin{frontmatter}
\title{A sequence approach to solve the Burgers' equation}
\author[mymainaddress]{S.A. Reijers}
\address[mymainaddress]{Physics of Fluids Group, Faculty of Science and Technology, MESA+Institute, University of Twente, P.O. Box 217,
7500 AE Enschede, The Netherlands}

\begin{abstract}
The Burgers' equation is a one-dimensional momentum equation for a Newtonian fluid. The Cole-Hopf transformation solves the equation for a given initial and boundary condition. However, in most cases the resulting integral equation can only be solved numerically. In this work a new semi-analytic solving method is introduced for analytic and bounded series solutions of the Burgers' equation. It is demonstrated that a sequence transformation can split the non-linear Burgers' equation into a sequence of linear diffusion equations. Each consecutive sequence element can be solved recursively using the Green's function method. The general solution to the Burgers' equation can therefore be written as a recursive integral equation for any initial and boundary condition. For a complex exponential function as initial condition we derive a new analytic solution of the Burgers' equation in terms of the Bell polynomials. The new solution converges absolutely and uniformly and matches a numerical solution with arbitrary precision. The presented semi-analytic solving method can be generalized to a larger class of non-linear partial differential equations which we leave for future work.
\end{abstract}

\begin{keyword}
Sequence transformation\sep Green's function method\sep nonlinear differential equations\sep Burgers' equation\sep Bell polynomials
\end{keyword}

\end{frontmatter}

\section{Introduction}
The Burgers' equation is a momentum equation for a one-dimensional viscous fluid \cite{Bonkile2018}. The equation describes a competition in the velocity field between non-linear convection and linear diffusion. It has an important historical significance in modelling turbulence, plays an important role in non-linear acoustics and has applications in traffic flows \cite{BURGERS1948171,murray_1973,sugimoto_1991,Rudenko2016, Musha_1976}.

The Cole-Hopf transformation converts the Burgers' equation into a linear diffusion equation which can be solved by the Green's function method for any initial and boundary condition \cite{Hopf1950,Cole2016OnAQ}. For some initial data the resulting integral equation can be solved explicitly in a closed-form solution or series solution \cite{RODIN1970401,BentonPlatzman}. For many initial and boundary conditions however explicit solutions cannot be obtained and the integral equation has to be solved numerically. As a result, the capability to analytically understand the non-linear dynamics is limited and other semi-analytic solving methods become important. 

In this work we solve the Burgers' equation using a sequence approach. We show that the Burgers' equation can be transformed into a sequence of linear diffusion equations for holomorphic and bounded series solutions. Each consecutive sequence element can be solved recursively using the Green's function method. The general solution to any initial value problem of the Burgers' equation can therefore be written as a recursive integral equation. This novel semi-analytic method is termed the sequence transformation method. For some initial and boundary conditions this recursive integral equation can be solved explicitly and written in a closed-form series. We show that we are able to obtain a closed-form series solution for a specific initial condition for which the Cole-Hopf integral solution holds no known closed-form solution.

Other semi-analytical methods such as the homotopy analysis method, Adomian decomposition method, differential transform method and variational iteration method have been successfully used in the past to solve the Burgers' equation \cite{SRIVASTAVA2014533,LIN20103082,BIAZAR20091394,ABBASBANDY20051265,AMINIKHAH20135979,PRAKASH2015314,BASTO2007927,ABDOU2005245}. The approach in these methods resembles our approach, in the sense that the solution is an infinite series whose terms can be calculated term-wise. However our approach does not require a deformation parameter, an auxiliary convergence parameter, predefined expansion polynomials, an expansion parameter, explicit Taylor series or Lagrange multipliers. Instead, the sequence of linear diffusion equations follows naturally from a sequence transformation applied to the Burgers' equation. Any initial condition, boundary condition or source term can be straightforwardly integrated using the Green's function method.

In \S \ref{sec:preliminaries} a Banach space is introduced including a set of sequence operations. The sequence transformation method is presented in \S \ref{sec:transformation}. In \S \ref{sec:burgerseq} the complex valued Burgers' differential equation is given. A new and complete analytic solution to the Burgers' equation is presented in \S \ref{sec:initialvalue} written in terms of the Bell polynomials for a complex exponential function as initial condition. Finally in \S \ref{sec:discussionandconclusion} we conclude and discuss possible future directions.

\section{Preliminaries\label{sec:preliminaries}}
In \S \ref{sub:banach} we introduce a Banach space for a class of well behaved and physically relevant functions. In the remainder of this section we discuss four important operations in the Banach space that are used in this work. First in \S\ref{sub:diffandint} we show that differentiation and integration for a series in this space can be performed term-wise. Finally in \S\ref{sub:product} the Cauchy product is introduced for two series.

\subsection{Banach space and series\label{sub:banach}}
Let $\Omega = \Omega_z \times \Omega_w \subseteq \mathbb{C}^2$ be an open subset and let $\mathcal{F}(\Omega)$ be a Banach space of all bounded holomorphic functions of two variables $f(z,w): \Omega \to \mathbb{C}$ with the supremum norm $||f(z,w)|| := \sup_{ (z,w) \in \Omega} |f(z,w)|$. We define the sequence space $\mathcal{F}_s(\Omega)$ as the set of all possible infinite sequences $\{f_n(z,w)\}$ with $f_n(z,w)\in\mathcal{F}(\Omega)$ and $n\in\mathbb{N}$ whose series converges by the Weierstrass M-test. That is, a sequence $\{f_n(z,w)\}$ is only member of $\mathcal{F}_s(\Omega)$ when $||f_n(z,w)|| \leq M_n$ and $\sum_{n=1}^\infty M_n < \infty$. As a result, any sequence $\{f_n(z,w)\}\in\mathcal{F}_s(\Omega)$ has a series that converges absolutely and uniformly on $\Omega$ to a member in $\mathcal{F}(\Omega)$, since the Banach space is complete. Our notation for a series is
\begin{equation}
f(z,w) := \sum \{f_n(z,w)\} := \sum_{n} f_n(z,w),
\label{eq:sequenceseriesdef}
\end{equation}
where the summation over a sequence results in a series.

\subsection{The differentiation and integration of a series\label{sub:diffandint}}
A series $f(z,w)\in\mathcal{F}(\Omega)$ of the sequence $\{f_n(z,w)\} \in \mathcal{F}_s(\Omega)$ is infinitely differentiable term-wise with respect to $z$ or $w$, since $f(z,w)$ converges uniformly in $\Omega$ and is analytic in both $z$ and $w$. Furthermore integration of the series $f(z,w)$ along any continuous path in $\Omega$ with respect to the variable $z$ or $w$ can be performed term-wise, since $f(z,w)$ converges uniformly in $\Omega$ and is continuous in both $z$ and $w$.

\subsection{The product of two series\label{sub:product}}
Two series $f(z,w)\in\mathcal{F}(\Omega)$ of the sequence $\{f_n(z,w)\} \in \mathcal{F}_s(\Omega)$ and $g(z,w)\in\mathcal{F}(\Omega)$ of the sequence $\{g_n(z,t)\} \in \mathcal{F}_s(\Omega)$ may be multiplied together using the Cauchy product as 
\begin{equation}
f(z,w)g(z,w) := \left(\sum_{n} f_n(z,w)\right)\left(\sum_{n} g_n(z,w)\right) = \sum_{n} \left(\sum_{m=1}^{n-1} f_m(z,w)g_{n-m}(z,w)\right).
\end{equation}
The functions $f(z,w)$, $g(z,w)$ and the product $f(z,w)g(z,w)$ are absolutely and uniformly convergent in $\Omega$. Therefore we have that the product $f(z,w)g(z,w)$ is a member of $\mathcal{F}(\Omega)$, since the Banach space is complete.

\section{The sequence transformation method\label{sec:transformation}}
Let us define a transformation $\mathcal{T}: \mathcal{F}_s(\Omega) \to \mathcal{F}(\Omega)$ for a sequence $\{f_n(z,w)\} \in \mathcal{F}_s(\Omega)$ to a tagged series $\hat{f}(z,w,s)\in\mathcal{F}(\Omega)$ by
\begin{align}
\mathcal{T}\bigg[\{f_n(z,w)\}\bigg] :=  \sum_{n} f_n(z,w)\exp(i n s),
\label{eq:transformationforward}
\end{align}
where $i$ is the imaginary unit and $s\in[-\pi,\pi]$. The tagged series 
\begin{equation}
\hat{f}(z,w,s) = \sum_{n} f_n(z,w)\exp(i n s)
\end{equation}
converges absolutely and uniform on $\Omega\times[-\pi,\pi]$, because
\begin{equation*}
||f_n(z,w)\exp(i n s)|| \leq ||f_n(z,w)|| \leq M_n
\end{equation*}
for all $(z,w,s)\in\Omega\times[-\pi,\pi]$. The inverse transformation $\mathcal{T}^{-1}: \mathcal{F}(\Omega)\to \mathcal{F}_s(\Omega)$ from a tagged series $\hat{f}(z,w,s)\in\mathcal{F}(\Omega)$ to a sequence $\{f_m(z,w)\} \in \mathcal{F}_s(\Omega)$ is defined as
\begin{equation}
\mathcal{T}^{-1}\bigg[\hat{f}(z,w,s)\bigg] = \left\{\frac{1}{2\pi}\int_{-\pi}^{\pi} \hat{f}(z,w,s)\exp(-i m s)\,\mathrm{d}s\right\} = \{f_m(z,w)\},
\end{equation}
where $m\in\mathbb{N}$ and we are allowed to integrate term-by-term, since the integrand converges uniformly on $\Omega\times[-\pi,\pi]$. The identity transform is given by
\begin{align*}
\mathcal{T}^{-1}\mathcal{T}\bigg[\{f_n(z,w)\}\bigg] &= \left\{\frac{1}{2\pi}\int_{-\pi}^{\pi} \left(\sum_{n} f_n(z,w)\exp(i n s)\right)\exp(-i m s)\,\mathrm{d}s\right\} \\&= \left\{\sum_{n}\frac{1}{2\pi}\int_{-\pi}^{\pi} f_n(z,w)\exp(i (n-m) s)\,\mathrm{d}s\right\} \\&= \{f_m(z,w)\}.
\end{align*}

The remarkable feature of this transformation is that every sequence in $\mathcal{F}_s(\Omega)$ has a tagged series for which the inverse transformation allows to obtain the original sequence. We can use this property to split a non-linear differential equation into a sequence of linear differential equations which we show in the next section for the Burgers' equation.

\section{The Burgers' equation\label{sec:burgerseq}}
Let $f(z,w)\in\mathcal{F}(\Omega)$ be a series of the sequence $\{f_n(z,w)\} \in \mathcal{F}_s(\Omega)$ and be a solution of the Burgers' differential operator $\mathcal{A}: \mathcal{F}(\Omega)\to\mathcal{F}(\Omega)$ given by the zero-map
\begin{align}
\mathcal{A}\bigg[f(z,w)\bigg] &= \frac{\partial f(z,w)}{\partial w} - \nu \frac{\partial^2 f(z,w)}{\partial z^2} + f(z,w)\frac{\partial f(z,w)}{\partial z}=0,\label{eq:burgers}
\end{align}
where $\nu\in\mathbb{R}>0$ is the viscosity.

In \S \ref{sub:sequencediffeq} we apply the sequence transformation to the Burgers' equation and show that the non-linear differential equation can be decomposed in a sequence of linear diffusion equations for all holomorphic and bounded series solutions. Finally in \S \ref{sub:greensfunction} the Green's function method is introduced to solve the resulting sequence of linear diffusion equations for an initial value problem and we present the general solution to the Burgers' equation.
\subsection{The sequence differential equation\label{sub:sequencediffeq}}
Instead of solving (\ref{eq:burgers}) directly for a function $f(z,w)$, we first decompose the non-linear equation into a sequence of linear diffusion equations by making use of the sequence transformation, see \S \ref{sec:transformation}. We start by applying the differential operator $\mathcal{A}$ to a tagged series $\hat{f}(z,w,s)$ and calculate the inverse sequence transformation $\mathcal{T}^{-1}$ by

\begin{align}
&\mathcal{T}^{-1}\bigg[\mathcal{A}\bigg[\hat{f}(z,w,s)\bigg]\bigg] = \mathcal{T}^{-1}\bigg[\frac{\partial \hat{f}(z,w,s)}{\partial w} - \nu \frac{\partial^2 \hat{f}(z,w,s)}{\partial z^2} + \hat{f}(z,w,s)\frac{\partial \hat{f}(z,w,s)}{\partial z}\bigg]\nonumber\\
&= \mathcal{T}^{-1}\bigg[\sum_{n} \bigg( \frac{\partial f_n(z,w)\exp(i n s)}{\partial w}  - \nu \frac{\partial^2 f_n(z,w)\exp(i n s)}{\partial z^2} \nonumber \\&+ \bigg(\sum_{l=1}^{n-1} f_l(z,w)\frac{\partial f_{n-l}(z,w)}{\partial z}\bigg)\exp(i n s)\bigg)\bigg]\nonumber\\
&= \left\{\frac{\partial f_m(z,w)}{\partial w} - \nu \frac{\partial^2 f_m(z,w)}{\partial z^2} + \sum_{l=1}^{m-1} f_l(z,w)\frac{\partial f_{m-l}(z,w)}{\partial z}\right\}\label{eq:recursiveseq}\\
&= \{0\}\nonumber,
\end{align}
where $\{0\}$ is the null sequence, $m\in\mathbb{N}$, the product $\hat{f}(z,w,s)\frac{\partial \hat{f}(z,w,s)}{\partial z}$ is written out using the Cauchy product (see section \ref{sub:product}) and we have used that $\hat{f}(z,w,s)$ is a holomorphic function that is absolutely and uniformly convergent. As a result, we see that each sequence element $f_m(z,w)$ is a solution of an in-homogeneous linear diffusion equation. 

In the next section we use the Green's function method to solve the sequence of linear diffusion equations.
\subsection{Solving the sequence differential equation in the complex plane\label{sub:greensfunction}}
In order to solve the series of in-homogenous linear differential equations (\ref{eq:recursiveseq}) using a Green's function, we must first define two open paths in the complex plane along which we will integrate the solution. Let $\gamma_z(t): [0,1] \to \Omega_z$ and $\gamma_w(t) : [0,1] \to \Omega_w$ be two continuous and piecewise differentiable open paths, where $\gamma_z(0) = z_a$, $\gamma_z(1) = z_b$, $\gamma_w(0) = w_a$ and $\gamma_w(1) = w_b$. We solve (\ref{eq:recursiveseq}) by multiplying the differential equation with $G(z,w;z_0,w_0): \Omega \times \Omega \to \mathbb{C}$, the Green's function, and integrate over the paths $\gamma_w$ and $\gamma_z$ to obtain
\begin{align}
\int_{\gamma_w} \int_{\gamma_z} G(z,w;z_0,w_0)\bigg( \frac{\partial f_m(z_0,w_0)}{\partial w_0} - \nu \frac{\partial^2 f_m(z_0,w_0)}{\partial z_0^2} \nonumber\\+ 
\sum_{l=1}^{m-1} f_l(z_0,w_0)\frac{\partial f_{m-l}(z_0,w_0)}{\partial z_0}\bigg)\,\mathrm{d}z_0\,\mathrm{d}w_0 = 0.
\end{align}
By using partial integration we can deduce that
\begin{align}
&\int_{\gamma_w} \int_{\gamma_z} f_m(z_0,w_0)\bigg( -\frac{\partial G(z,w;z_0,w_0)}{\partial w_0} - \nu \frac{\partial^2 G(z,w;z_0,w_0)}{\partial z_0^2}\bigg)\,\mathrm{d}z_0\,\mathrm{d}w_0 = \nonumber\\
&-\int_{\gamma_w} \int_{\gamma_z}  G(z,w;z_0,w_0)\bigg(\sum_{l=1}^{m-1} f_l(z_0,w_0)\frac{\partial f_{m-l}(z_0,w_0)}{\partial z_0}\bigg)\,\mathrm{d}z_0\,\mathrm{d}w_0 \nonumber +\\
&\nu\int_{\gamma_w} \left[G(z,w;z_0,w_0)\frac{\partial f_m(z_0,w_0)}{\partial z_0}-\frac{\partial G(z,w;z_0,w_0)}{\partial z_0}f_m(z_0,w_0)\right]_{z_0=z_a}^{z_0=z_b} \,\mathrm{d}w_0 \nonumber \\
&-\int_{\gamma_z} \left[G(z,w;z_0,w_0)f_m(z_0,w_0)\right]_{w_0=w_a}^{w_0=w_b} \,\mathrm{d}z_0. \label{eq:solpartialint}
\end{align}
If the Green's function $G(z,w;z_0,w_0)$ satisfies 
\begin{align}
-\frac{ \partial G(z,w;z_0,w_0)}{\partial w_0} - \nu\frac{\partial^2 G(z,w;z_0,w_0)}{\partial z_0^2} = \delta(z-z_0)\delta(w-w_0),
\end{align}
where $\delta$ is the complex Dirac delta function, we can write rewrite (\ref{eq:solpartialint}) as
\begin{align}
&f_m(z,w) = -\int_{\gamma_w} \int_{\gamma_z}  G(z,w;z_0,w_0)\bigg(\sum_{l=1}^{m-1} f_l(z_0,w_0)\frac{\partial f_{m-l}(z_0,w_0)}{\partial z_0}\bigg)\,\mathrm{d}z_0\,\mathrm{d}w_0 \nonumber +\\
&\nu\int_{\gamma_w} \left[G(z,w;z_0,w_0)\frac{\partial f_m(z_0,w_0)}{\partial z_0}-\frac{\partial G(z,w;z_0,w_0)}{\partial z_0}f_m(z_0,w_0)\right]_{z_0=z_a}^{z_0=z_b} \,\mathrm{d}w_0 \nonumber \\
&-\int_{\gamma_z} \left[G(z,w;z_0,w_0)f_m(z_0,w_0)\right]_{w_0=w_a}^{w_0=w_b} \,\mathrm{d}z_0, \label{eq:generalsol}
\end{align}
if and only if the paths $\gamma_w(t)$ and $\gamma_z(t)$ cross the Dirac delta function at $z=z_0$ and $w=w_0$. We note that when the integrand is holomorphic on $\Omega$ the exact paths $\gamma_w(t)$ and $\gamma_z(t)$ do not play a role since any closed path in $\Omega$ integrates to zero for holomorphic integrands, see Morera's theorem \cite{Arfken1985}.

In order to integrate all boundary and initial conditions only once in the general solution to the Burgers' equation $f(z,w)$, we set
\begin{align}
&f_1(z,w) = \nu\int_{\gamma_w} \left[G(z,w;z_0,w_0)\frac{\partial f(z_0,w_0)}{\partial z_0}-\frac{\partial G(z,w;z_0,w_0)}{\partial z_0}f(z_0,w_0)\right]_{z_0=z_a}^{z_0=z_b} \,\mathrm{d}w_0 \nonumber \\
&-\int_{\gamma_z} \left[G(z,w;z_0,w_0)f(z_0,w_0)\right]_{w_0=w_a}^{w_0=w_b} \,\mathrm{d}z_0, \nonumber \\
&f_{m\ge 2}(z,w) = -\int_{\gamma_w} \int_{\gamma_z}  G(z,w;z_0,w_0)\bigg(\sum_{l=1}^{m-1} f_l(z_0,w_0)\frac{\partial f_{m-l}(z_0,w_0)}{\partial z_0}\bigg)\,\mathrm{d}z_0\,\mathrm{d}w_0. \label{eq:generalfm}
\end{align}
We note that this choice is not unique and other ways are possible to integrate the boundary and initial conditions in the solution which is outside the scope of this work. The general solution of the Burgers' equation (\ref{eq:burgers}) for any initial and boundary condition can now be written as a recursive integral equation
\begin{align}
&f(z,w) = \sum_{m=1}^\infty f_m(z,w) = \nonumber\\ &-\sum_{m=2}^\infty\int_{\gamma_w} \int_{\gamma_z}  G(z,w;z_0,w_0)\bigg(\sum_{l=1}^{m-1} f_l(z_0,w_0)\frac{\partial f_{m-l}(z_0,w_0)}{\partial z_0}\bigg)\,\mathrm{d}z_0\,\mathrm{d}w_0 \nonumber +\\
&\nu\int_{\gamma_w} \left[G(z,w;z_0,w_0)\frac{\partial f(z_0,w_0)}{\partial z_0}-\frac{\partial G(z,w;z_0,w_0)}{\partial z_0}f(z_0,w_0)\right]_{z_0=z_a}^{z_0=z_b} \,\mathrm{d}w_0 \nonumber \\
&-\int_{\gamma_z} \left[G(z,w;z_0,w_0)f(z_0,w_0)\right]_{w_0=w_a}^{w_0=w_b} \,\mathrm{d}z_0. \label{eq:generalsolb}
\end{align}
\section{An initial value problem for the Burgers' equation\label{sec:initialvalue}}
Consider the following initial value problem for the Burgers' equation
\begin{gather}
\frac{\partial f(z,w)}{\partial w} - \nu \frac{\partial^2 f(z,w)}{\partial z^2} + f(z,w)\frac{\partial f(z,w)}{\partial z}=0\nonumber,\\
f(z,0) = \exp(i z) = \cos(z) + i \sin(z)\label{eq:initialvalueproblem}
 \end{gather}
where $\nu \in \mathbb{R}>0$, $\Omega_w=\{\Re(w)\ge 0,\, w\in \mathbb{C}\}$ and $\Omega_z=\{\Im(z)\ge 0,\, z\in \mathbb{C}\}$. To solve the differential equation, we start in \S \ref{sub:solving} with the ansatz that the solution is a series that can be solved using the sequence transformation method. Later in \S \ref{sub:convergence} we proof that the ansatz is correct and we show that the series indeed converges absolutely and uniform in a domain $\Omega$. Finally in \S \ref{sub:abserror} we give the absolute error of the series solution and compare the solution convergence properties to a numerical solution using the Cole-Hopf transformation.

\subsection{Solving the initial value problem\label{sub:solving}}
We assume that the solution to the initial value problem (\ref{eq:initialvalueproblem}) is a series $f(z,w)\in\mathcal{F}(\Omega)$ of the sequence $\{f_n(x,t)\} \in \mathcal{F}_s(\Omega)$. Following (\ref{eq:generalfm}), the solution for the sequence $\{f_m(z,w)\}$ is
\begin{align*}
f_1(z,w) &= \int_{\gamma_z} G(z,w;z_0,w_0) \exp(i z_0)\,\mathrm{d}z_0,\\
f_{m \ge 2}(z,w) &= -\int_{\gamma_w}\int_{\gamma_z} G(z,w;z_0,w_0) \left(\sum_{l=1}^{m-1} f_l(z_0,w_0)\frac{\partial f_{m-l}(z_0,w_0)}{\partial z_0}\right)\,\mathrm{d}z_0\,\mathrm{d}w_0,
\end{align*}
where $G(z,w;z_0,w_0)$ is the free-space Green's function
\begin{equation}
G(z,w;z_0,w_0) = \frac{1}{\sqrt{4\pi \nu (w-w_0)}}\exp\left(\frac{-(z-z_0)^2}{4\nu(w-w_0)}\right)H(w-w_0),
\end{equation}
$i$ is the imaginary unit, $H(w-w_0)$ is the Heaviside theta function, the path $\gamma_w(t)$ connects $w_0=0$ to $w_0=w^+$ and the path $\gamma_z(t)$ connects $z_0=-\infty$ to $z_0=\infty$ \cite{Duffy2015}. We note that the integrand is holomorphic on the domain of interest and thus the exact paths do not matter as long they cross the Dirac delta function at $z=z_0$ and $w = w_0$. The first three calculated sequence elements of $\{f_m(z,w)\}$ are 
\begin{align*}
&\begin{dcases}
f_1(z,w) &= \exp(i z - w \nu),\\
f_2(z,w) &= \frac{i  \exp(-4 \nu  w+2 i z)\left(1-\exp(2 \nu  w)\right)}{2 \nu },\\
f_3(z,w) &=-\frac{\exp(-9 \nu  w+3 i z)\left(1-3 \exp(4 \nu  w)+2 \exp(6 \nu  w)\right)}{8 \nu ^2},\\
\cdots &=\cdots
\end{dcases}
\end{align*}
where we have omitted the Heaviside theta functions since $w\ge0^+$. The sequence can be captured in a closed from as
\begin{align}
f_m(z,w) =& \frac{i^{m-1}  \exp(-\nu  m^2 w+i m z) }{2^{m-1} \nu ^{m-1}(m-1)!}  \times\nonumber\\& \sum _{k=1}^m (-1)^{k-1} (k-1)! \text{B}_{m,k}\left(\mu_1(m,w),\dots,\mu_{m-k+1}(m,w)\right),
\label{eq:fundamentalseqsolutionintialvalueproblem}
\end{align}
where $\text{B}_{m,k}\left(\mu_1(m,w),\dots,\mu_{m-k+1}(m,w)\right)$ are the exponential Bell polynomials and
\begin{equation*}
\mu_l(m,w) := \exp (l \nu   w (m-l)).
\end{equation*}

\begin{figure}[t]
        \centering
                \includegraphics[width=\textwidth]{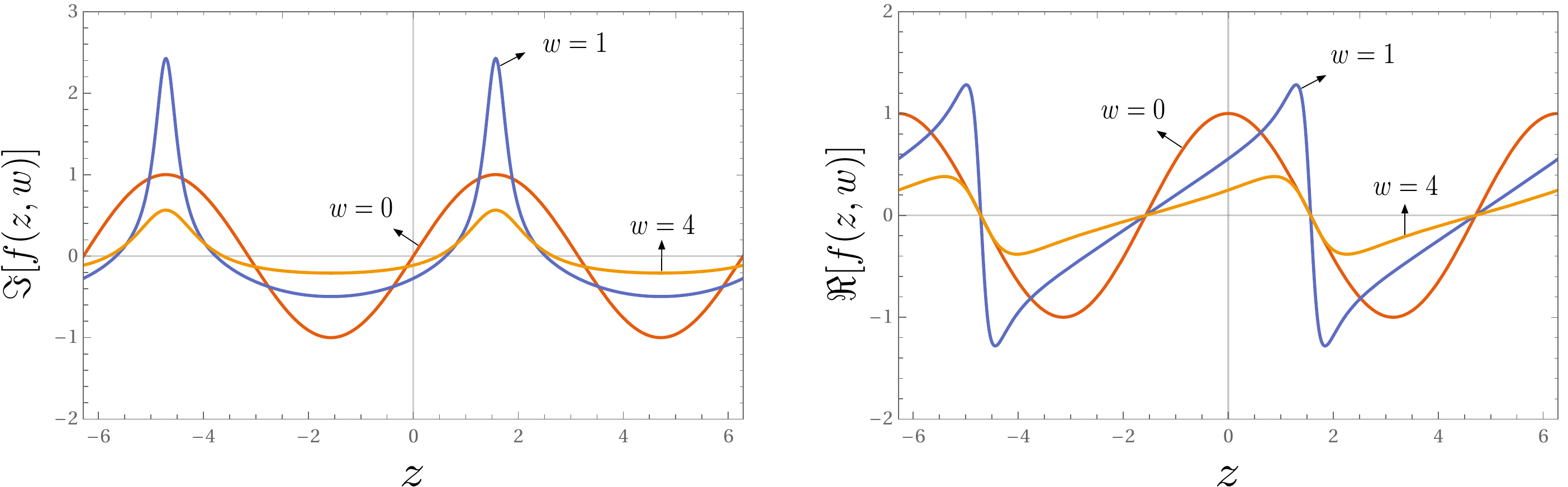}
                \caption{A plot of the complex-valued solution $f(z,w)$, given by (\ref{eq:fundamentalsolutionintialvalueproblem}), for $2\pi < z < 2 \pi$ at three different instances $w=0$, $w=1$ and $w=4$ with $\nu=0.3$. The imaginary part is plotted on the left and the real part is plotted on the right. The series solution has been cut-off at $m=30$.}
                \label{fig:solutionplot}
\end{figure}

Following (\ref{eq:generalsolb}), the full solution to our initial value problem is then given by
\begin{equation}
f(z,w) = \sum_{m=1}^\infty f_m(z,w),
\label{eq:fundamentalsolutionintialvalueproblem}
\end{equation}
where $(z,w)\in\Omega$. In figure \ref{fig:solutionplot} we show the solution $f(z,w)$ for $2\pi < z < 2\pi$ at different instances with $\nu=0.3$. The figure shows a competition between non-linear convection and linear diffusion in the velocity field. At $w=1$ we see the development of a shock wave in both the real and imaginary part of the solution which is dissipated by viscosity on later times.

\subsection{Proof of convergence\label{sub:convergence}}
In order to show that our ansatz was correct, we need to proof that (\ref{eq:fundamentalsolutionintialvalueproblem}) converges absolutely and uniform in $\Omega$. Instead, here we show that the solution converges on the real line and leave the general proof in the complex plane for future work. We introduce $u(x,t) : \Delta_x \times\Delta_t \to \mathbb{C}$ defined as $u(x,t) := f(x,t)$, where $\Delta_x = \{x\in\mathbb{R}\}$ and $\Delta_t = \{t\ge 0,\,t\in\mathbb{R}\}$. We start by
\begin{align}
|u_m(x,t)| &\leq \frac{\exp(-\nu  m^2 t)}{2^{m-1} \nu ^{m-1}(m-1)!} \bigg| \sum _{k=1}^m (-1)^{k-1} (k-1)! \text{B}_{m,k}\left(\mu_1(m,t),\dots,\mu_{m-k+1}(m,t)\right)\bigg|,\nonumber\\
&\leq\frac{\exp(-\nu  m^2 t)}{2^{m-1} \nu ^{m-1}(m-1)!} \sum _{k=1}^m \bigg| (-1)^{k-1} (k-1)! \text{B}_{m,k}\left(\mu_1(m,t),\dots,\mu_{m-k+1}(m,t)\right)\bigg|,\nonumber\\
&=\frac{\exp(-\nu  m^2 t)}{2^{m-1} \nu ^{m-1}(m-1)!} \sum _{k=1}^m (k-1)! \text{B}_{m,k}\left(\mu_1(m,t),\dots,\mu_{m-k+1}(m,t)\right),\nonumber\\
&\leq\frac{\exp(-\nu  m^2 t)\exp(\nu (m-1)m t)}{2^{m-1} \nu ^{m-1}(m-1)!} \sum _{k=1}^m  (k-1)! \text{B}_{m,k}\left(1,1,\dots,1\right),\nonumber\\
&=\frac{\exp(-\nu m t)}{2^{m-1} \nu ^{m-1}(m-1)!} \sum _{k=1}^m  (k-1)!{m\brace k},\label{eq:estimation}
\end{align}
where ${m\brace k}$ are the Stirling numbers of the second kind. Now by the ratio test we have that
\begin{align}
\lim_{m\to\infty} \frac{|u_{m+1}(x,t)| }{|u_m(x,t)| } &= \lim_{m\to\infty}\left(\frac{\frac{\exp(-\nu  (m+1) t)}{2^{m} \nu ^{m}(m)!} \sum _{k=1}^{m+1}  (k-1)!{m+1\brace k}}{\frac{\exp(-\nu  m t)}{2^{m-1} \nu ^{m-1}(m-1)!} \sum _{k=1}^m  (k-1)!{m\brace k}}\right)\nonumber\\
&= \lim_{m\to\infty}\left( \frac{\exp(-\nu t)}{2\nu m} \frac{\sum _{k=1}^{m+1}  (k-1)!{m+1\brace k}}{\sum _{k=1}^m  (k-1)!{m\brace k}}\right),\nonumber\\
&= \frac{\exp(-\nu t)}{2\nu}\lim_{m\to\infty}\left( \frac{1}{m} \frac{\sum _{k=1}^{m+1}  (k-1)!{m+1\brace k}}{\sum _{k=1}^m  (k-1)!{m\brace k}}\right),\nonumber\\
&= \frac{r \exp(-\nu t)}{2\nu},\nonumber\\
 &< 1,\label{eq:ratiotest}
\end{align}
where $r = \lim_{m\to\infty}\left( \frac{1}{m} \frac{\sum _{k=1}^{m+1}  (k-1)!{m+1\brace k}}{\sum _{k=1}^m  (k-1)!{m\brace k}}\right)$ and can be calculated numerically as $r\approx1.4427$. The convergence condition for the Weierstrass M-test is fulfilled, see section \ref{sub:banach}, when $\nu > r/2$ for all $t\ge 0$.
\begin{figure}[t]
        \centering
                \includegraphics[width=0.8\textwidth]{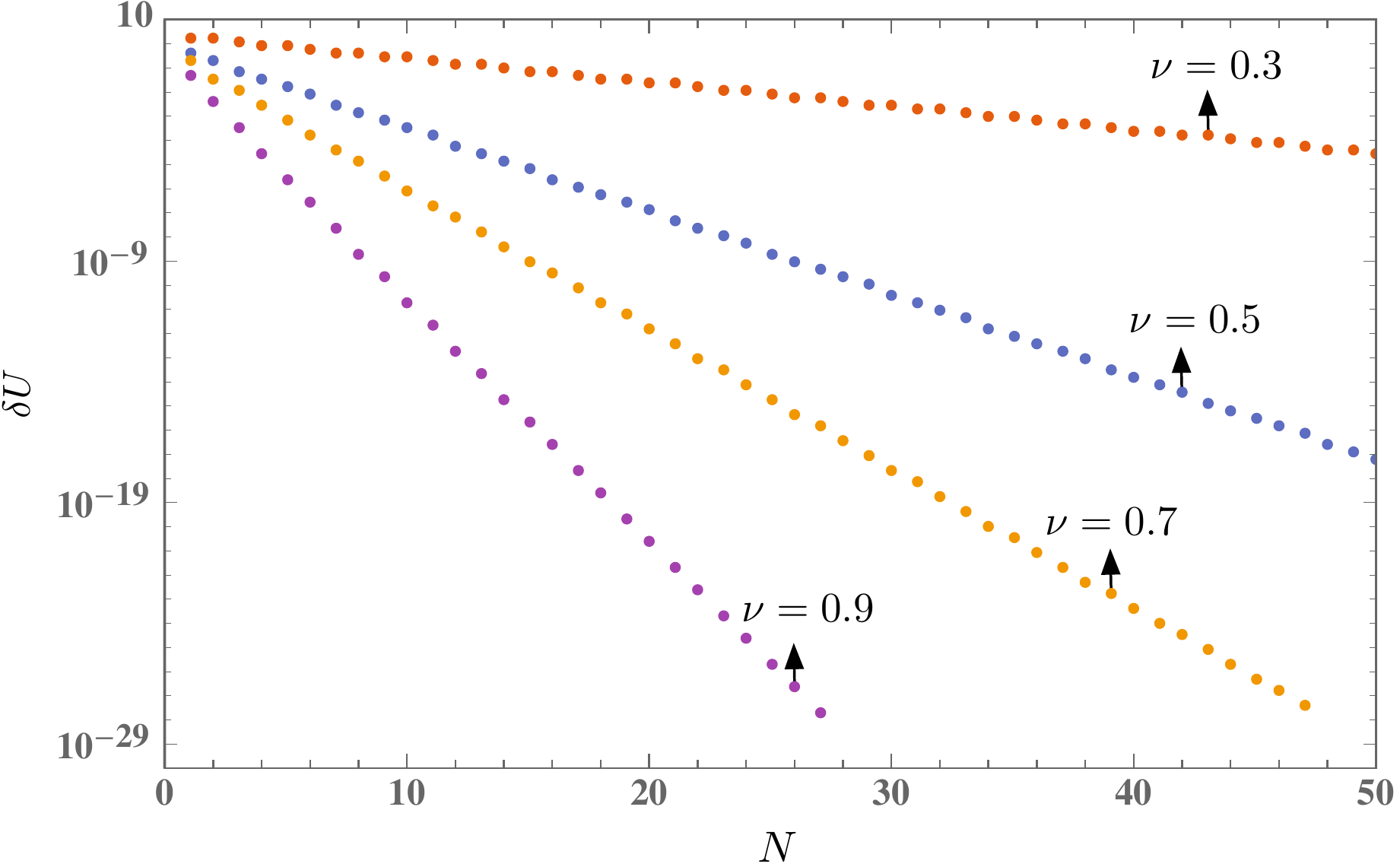}
                \caption{A log plot of the absolute error between a high-precision numerical solution ($30$-digits accurate) of (\ref{eq:colehopfresult}) and the partial series solution (\ref{eq:fundamentalsolutioncutoff}) as function of $N$ for different viscosities in $\Lambda= [-2\pi,2\pi]\times[0,3]$, see (\ref{eq:abserror}).}
                \label{fig:convergence}
\end{figure}
\begin{figure}[t]
        \centering
                \includegraphics[width=0.8\textwidth]{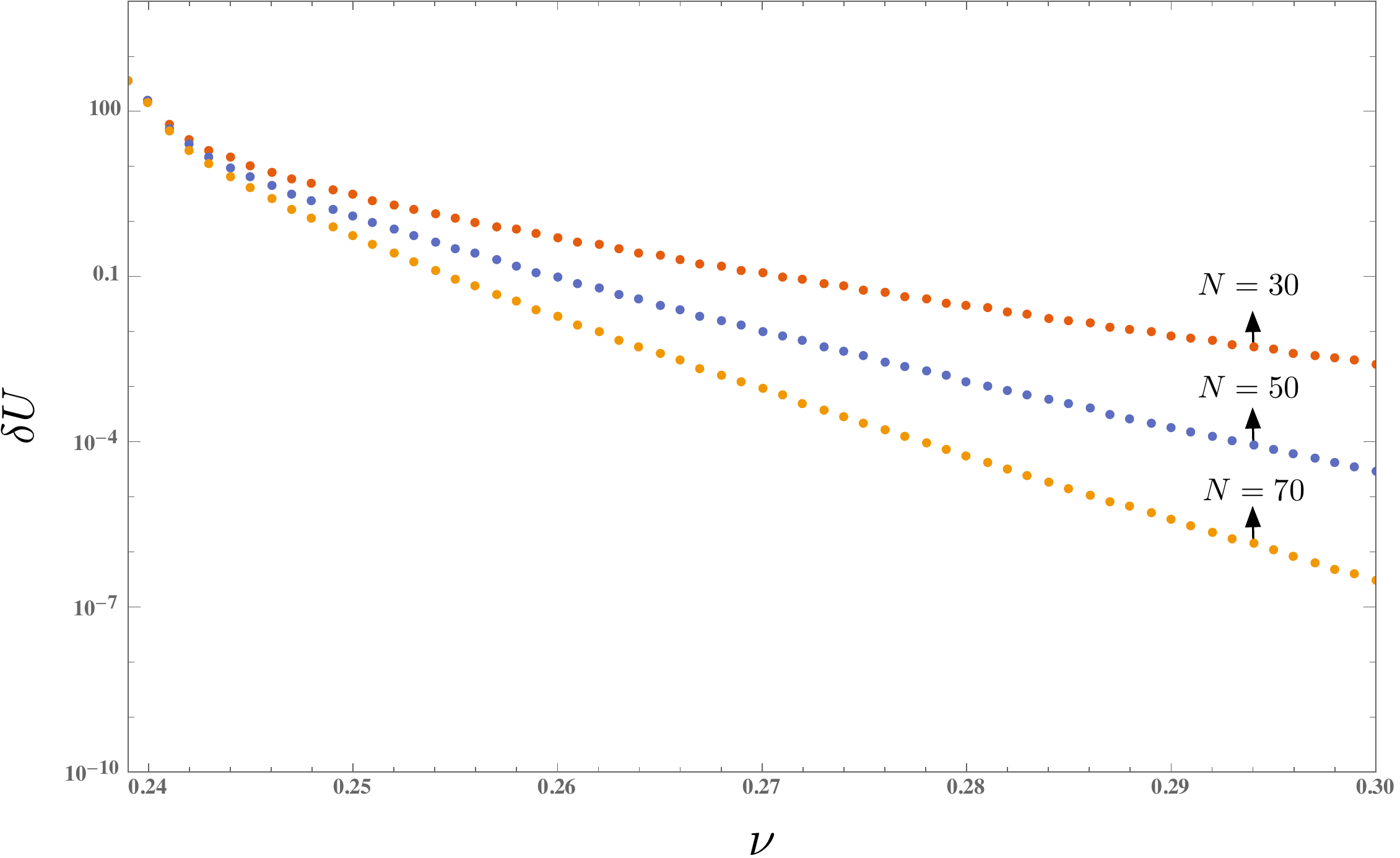}
                \caption{A log plot of the absolute error between a high-precision numerical solution ($30$-digits accurate) of (\ref{eq:colehopfresult}) and the partial series solution (\ref{eq:fundamentalsolutioncutoff}) as function of the viscosity $\nu$ for different $N$ in $\Lambda= [-2\pi,2\pi]\times[0,3]$, see (\ref{eq:abserror}).}
                \label{fig:singularity}
\end{figure}
We note that we do not have a more narrow upper bound for the Bell polynomials in (\ref{eq:estimation}), therefore the series may in fact be convergent for values $\nu \le r/2$ as figure \ref{fig:solutionplot} shows. 

\subsection{The absolute error\label{sub:abserror}}
The series solution (\ref{eq:fundamentalsolutionintialvalueproblem}) of the initial value problem (\ref{eq:initialvalueproblem}) converges absolutely and uniform in $\Delta= (-\infty,\infty)\times[0,\infty)$ for $\nu > r/2$. This means that we can define an absolute error as
\begin{align}
\delta U = \left|U(x,t)-U_N(x,t)\right|,\label{eq:abserror}
\end{align}
where $(x,t) = \sup_{ (x,t) \in  \Lambda}|U(x,t)|$, $\Lambda\subseteq \Delta$, $U(x,t)$ is the exact solution of the initial value problem and $U_N(x,t)$ is the partial series solution given by 
\begin{equation}
U_N(x,t) = \sum_{m=1}^N u_m(x,t),
\label{eq:fundamentalsolutioncutoff}
\end{equation}
$N \in\mathbb{N}$. The absolute error (\ref{eq:abserror}) goes to zero as $N\to\infty$ and can be used to show how fast the partial series solution converges to the exact solution in $\Lambda$ for increasing values of $N$. 

The exact solution with infinite precision is not known, but we can use the absolute error (\ref{eq:abserror}) to show how fast the partial series solution (\ref{eq:fundamentalsolutioncutoff}) converges to the numerical solution of the initial value problem (\ref{eq:initialvalueproblem}) given by the Cole-Hopf transformation \cite{Hopf1950}
\begin{equation}
U(x,t) = \frac{\int_{-\infty}^{\infty} \frac{(x-x_0)}{t} \exp\left(\frac{-(x-x_0)^2}{4\nu t}-\frac{1}{2\nu}\int_0^{x_0} \exp(i x')\,\mathrm{d}x' \right)\,\mathrm{d}x_0}{\int_{-\infty}^{\infty} \exp\left(\frac{-(x-x_0)^2}{4\nu t}-\frac{1}{2\nu}\int_0^{x_0} \exp(i x')\,\mathrm{d}x'\right)\,\mathrm{d}x_0}.
\label{eq:colehopfresult}
\end{equation}
This integral does not have a known closed-form expression, but we can numerically estimate this integral with high-precision ($30$ digits) using Mathematica \cite{WolframMathematica}. 

Figure \ref{fig:convergence} shows a log plot of the absolute error between a high-precision numerical solution of (\ref{eq:colehopfresult}) and the partial series solution (\ref{eq:fundamentalsolutioncutoff}) as function of $N$ for different viscosities in $\Lambda= [-2\pi,2\pi]\times[0,3]$. The figure shows that the absolute error decreases exponentially with increasing $N$. Interestingly, the figure also shows that the rate of convergence decreases for decreasing $\nu$. From a physical point of view we expect this to happen, since in the limit $\nu\to 0$ the Burgers' equation has a finite time singularity where the solution blows up. Figure \ref{fig:singularity} shows a log plot of the absolute error as function of the viscosity for different $N$. This figure shows that indeed the viscosity plays an important factor in the convergence of the series solution. A singularity occurs at $\nu\approx 0.239$ for all plotted values $N$ where absolute error starts to grows exponentially. In this case the solution has a finite blow-up time where velocity gradient becomes infinitely steep at some point in $\Lambda= [-2\pi,2\pi]\times[0,3]$. Furthermore, the numerical estimation of (\ref{eq:colehopfresult}) also starts to produce convergence errors around this value for the viscosity. We discussed in \S \ref{sub:convergence} that we need a more narrow upper bound for the Bell polynomials in (\ref{eq:estimation}) in order to proof that the series is indeed convergent for values $\nu \le r/2$, which is outside the scope of this work. 
\section{Discussion \& Conclusion\label{sec:discussionandconclusion}}
In this work we introduced a novel semi-analytic method to solve the Burgers' equation. We have shown that for any holomorphic and bounded series solution, the Burgers' equation can be decomposed into a sequence of linear diffusion equations by means of a sequence transformation. The solution for each individual sequence element can be written as a recursive integral equation using the Green's function method. In some cases, this recursive integral can be explicitly solved and a closed-form expression for the sequence can be obtained.

Using the sequence transformation method, we solved an initial value problem of the Burgers' equation in the complex plane. We showed that, for a complex exponential function as initial condition, the solution can be written as an infinite series using the Bell polynomials. The series solution is analytically complete, bounded and converges absolutely when the viscosity is larger than a threshold value. Furthermore we showed that the series solution converges exponentially to a high-precision numerical solution using the Cole-Hopf transformation for specific values of the viscosity. 

In many physical systems one expects the solution of a non-linear differential equation to be analytic or at-least piece-wise analytic in the domain of interest. Therefore, we anticipate the sequence transformation method to be applicable far beyond the Burgers' equation. It would be very interesting to generalize the sequence transformation method to a broader set of non-linear differential equations that admit analytic solutions, e.g. wave equations such as the Sine Gordon equation. Furthermore, it would also be interesting to investigate the possibility of capturing the initial conditions and boundary conditions in a set of basis functions that form an orthonormal set with respect to the domain of interest. This way, the recursive integrals that integrate over this domain simply vanish or can be drastically simplified. Therefore, the general solution could possibly be written without the need of any integrals. This direction could open the doors to more efficient numerical solvers for non-linear differential equations without the need of any spatial or time discretization, which we leave for future work.

\section{Acknowledgements}
We are grateful to the OEIS Foundation Inc. for maintaining an extremely useful online encyclopedia of integer sequences, \url{http://oeis.org}. This work is part of an Industrial Partnership Programme of the Netherlands Organization for Scientific Research (NWO). This research programme is co-financed by ASML.

\bibliography{mybibfile}

\end{document}